\documentclass[12pt]{article}
\usepackage{xy}
\usepackage[latin1]{inputenc}
\usepackage[english]{babel}
\usepackage{amssymb,amsmath,amscd,amstext}
\usepackage{amsfonts}
\usepackage{amsthm}
\usepackage[dvips]{graphicx}
\usepackage{epsfig}
\usepackage{palatino}
\input{xy}
\xyoption{all}

\newcounter{eno}
\newcounter{ch}
\newcounter{sec}
\newcounter{no}
\def \sec#1 {\stepcounter{sec} \setcounter{eno}{1} \setcounter{no}{1}}
\def \prop#1 {\par{\bigskip \noindent \bf Proposition \arabic{sec}.\arabic{no}} \stepcounter{no} {\it #1}}
\def \cor#1 {\par{\bigskip \noindent \bf Corollary \arabic{sec}.\arabic{no}} \stepcounter{no} {\it #1}}
\def \lem#1 {\par{\bigskip \noindent \bf Lemma \arabic{sec}.\arabic{no}} \stepcounter{no} {\it #1}}
\def \thm#1 {\par{\bigskip \noindent \bf Theorem \arabic{sec}.\arabic{no}} \stepcounter{no} {\it #1}}
\def \conj#1 {\par{\bigskip \noindent \bf Conjecture \arabic{sec}.\arabic{no}} \stepcounter{no} {\it #1}}
\def \define#1 {\par{\bigskip \noindent \bf Definition \arabic{sec}.\arabic{no}} \stepcounter{no} {\it #1}\hspace{2mm}}

\newcommand {\Z} {\mathbb Z}
\newcommand {\Q} {\mathbb Q}
\newcommand {\ZG} {\mathbb {Z}[G]}
\newcommand {\QG} {\mathbb {Q}[G]}

\newcommand {\rk} {{\rm rk}_{\Z}}

\def \proof#1 {Proof:$\,\,\,$  #1   \newline ${}$ \hfill $\Box$
\,\,  }

\begin{document}

\begin{center}
{\huge Periodic Cohomology}

{W.H.Mannan}
\end{center}

\bigskip
{MSC16E05.}\hfill {Keywords:  projective resolution, periodic homology}

\bigskip
\begin{tiny}
\noindent  {\bf Abstract} We offer a direct proof of an elementary result
concerning cohomological periods.  As a corollary we show that
given a finitely generated stably free resolution of $\Z$ over a
finite group, two of its modules are free.
\end{tiny}

 \stepcounter{sec}

\sec{}
\section{Introduction}
The study of periodic cohomology over finite groups closely
relates to free group actions on spheres (see for example
\cite{Miln1}). This in part contributes to the importance of the
classification of finite groups with periodic cohomology (see
\cite{Suzu} and \cite{Zass}).

A much simpler classification is the classification of groups which
specifically have cohomological period 2. They must be cyclic.
Although well known, the proof in the literature (see \cite{Swan})
is disproportionately technical, involving the methods used in more
general classifications of periodic cohomology.
 This includes considering potential Sylow $p$-subgroups and finding normal
complements.  We offer a direct proof avoiding such technical
difficulties.

In \S2 we use this result to prove a result regarding finite stably
free resolutions of $\Z$ over finite groups:

\bigskip
\noindent {\bf Theorem A}  {\it Let $G$ be a finite group and
suppose we have a finitely generated stably free resolution over
$G$:
$$
 \cdots \stackrel
{d_3}\rightarrow S_2 \stackrel {d_2}\rightarrow S_1 \stackrel{d_1}
\rightarrow S_0 \rightarrow \Z
$$ Then $S_1$ and $S_2$ are free.}

\bigskip
The chain complex associated to the universal cover of a
cell complex may be regarded as an algebraic complex over the
fundamental group. In particular if the cell complex is finite and
2 dimensional one obtains a truncated free finite resolution of
$\Z$.  The classification of these homotopy types is of particular
interest due to its relation to Wall's D2 problem (see the
introduction to \cite{John1}).

Although the modules in these complexes are free, algebraic
surgery may leave some modules potentially only stably free. From
the point of view of classification it is useful to know that
these modules will still be free (in particular as the maps may
then be represented by matrices). Over finite groups Theorem A
does just that.

\sec{} \section{Cohomological Period}

We follow the characterization of cohomological period in
\cite{John1} by saying $n>0$ is a (cohomological) period of a
finite group $G$ if one of the following equivalent conditions
hold (see \cite{John1}, \S40, $\mathcal{P}_1(n)$,
$\mathcal{P}_3(n)$, $\mathcal{P}_4(n)$):
\newline\newline
$\mathcal{C}_1(n)$:\hspace{1mm}  $\mathcal{D}_{n+k}(\Z) \cong
\mathcal{D}_{k}(\Z)$ for all integers $k$.  ($\mathcal{D}_{n}$
denotes the $n^{\rm th}$ derived functor (see \cite{John1}, \S20) so
in particular $\widehat{H}^{an+k}(G;\,\Z) \cong
\widehat{H}^{bn+k}(G;\,\Z)$ for all $a,b,k \in \Z$).
\newline\newline
$\mathcal{C}_2(n)$:\hspace{1mm}  There exists an exact sequence of
the form:

$$\,\,\,\, 0 \to \Z \to P_{n-1}\to \cdots \to P_0 \to \Z \to 0$$
$\mathcal{C}_3(n)$:\hspace{1mm} $\widehat{H}^{n}(G;\,\Z) \cong \Z
/|G|$.  \hspace{2mm}(Note $\mathcal{C}_1(n)$ implies
$H_{n-1}(G;\,\Z) \cong \widehat{H}^{-n}(G;\,\Z)$

\noindent$ \cong \widehat{H}^{n}(G;\,\Z)$. In particular
$H_{n-1}(G; \, \Z) \cong \Z /|G|$ if $n$ is a period of $G$).

\thm{If 2 is a period of a finite group $G$ then $G$ is cyclic.}

\noindent Proof:\hspace{2mm} $G/ G' \cong H_1(G;\,\Z) \cong \Z
/|G|$ which has the same order as $G$ so $G \cong G/G'$.

${}$\hfill $\Box$ \,\,

\sec{} \section{Stably free resolutions}

Before using theorem 2.1 to prove Theorem A, we note the following
restriction on stably free modules over finite groups:

\prop{Any stably free module of finite $\ZG$- rank greater than one
is free.}

\noindent \proof{$\ZG \oplus \ZG$ is an Eichler lattice so the
result follows from \cite{John1}, theorem 15.1.}

\bigskip
\noindent Proof of Theorem A:\hspace{2mm}  {Suppose one of $S_1$ or
$S_2$ is not free and let $K \cong {\rm ker}(d_1)$.  We have an
exact sequence:
$$
0 \to K \rightarrow S_1 \stackrel{d_1}{\rightarrow}S_0  \to \Z \to
0 \eqno(1)
$$
Consideration of ranks and nullities implies $\rk(K) \equiv 1$ Mod
$|G|$.
\newline
${}$\hspace{2mm}$K$ is a submodule of $S_1$ and $d_2$ induces a
surjection $S_2 \to K$.  Hence $\rk(K) \leq \rk(S_1),\, \rk(S_2)$.
By proposition 3.1 we have $\rk(S_1) \leq |G|$ or $\rk(S_2) \leq
|G|$ so $\rk(K) = 1$.
\newline
${}$\hspace{2mm}Tensoring (1) with $\Q$ yields the exact sequence:
$$
0 \to K \otimes \Q \rightarrow \QG^a {\rightarrow} \QG^b  \to \Q
\to 0
$$
 By the 'Whitehead Trick', $K \otimes \Q \oplus \QG^b \cong
\Q\oplus \QG^a$  as $\QG$- modules. Consideration of dimension over
$\Q$ implies $a=b$ so cancellation gives $K \otimes \Q \cong \Q$. In
particular the $G$- action on $K \otimes \Q$ is trivial.  Hence the
$G$- action on $K$ is trivial and we have $K \cong \Z$. We therefore
have an exact sequence:
$$
0 \to \Z \rightarrow S_1 \stackrel{d_1}{\rightarrow}S_0  \to \Z
\to 0
$$
Hence $G$ satisfies $\mathcal{C}_2(2)$ and is cyclic by theorem
2.1.

Any finitely generated stably free module over a cyclic group is
free (see \cite{John1}, proposition 15.7). So in particular $S_1$
and $S_2$ are free as required. ${}$\hfill $\Box$ \,\,}

\bigskip
Note however that $S_0$ need not be free. For example if $G$ is the
quaternionic group $Q_{32}$ we have a finitely generated stably free
resolution where $S_0$ is not free (see \cite{John2}, \S 4).

\noindent email: \verb|wajid@mannan.info|

\bigskip
\noindent Address:

 \noindent   J8 Hicks Building\\
    Department of Pure Mathematics\\
    University of Sheffield\\
    Hounsfield Road\\
    Sheffield S3 7RH\\

\end{document}